\documentclass[11pt]{article}
\usepackage{amsfonts,amscd,amssymb, amsmath}
\newtheorem{definition}{Definition}[section]

\newtheorem{proposition}[definition]{Proposition}

\newtheorem{theorem}[definition]{Theorem}

\def\mf{\mathfrak}

\def\rawo\lonra{\longrightarrow}

\def\ot{\otimes}

\newcommand{\selabel}[1]{\label{se:#1}}

\allowdisplaybreaks[4]

\begin{document}
\title{More examples of invariance under twisting 
\thanks{Research partially supported by the CNCSIS project 
 ''Hopf algebras, cyclic homology and monoidal categories'', 
contract nr. 560/2009, CNCSIS code $ID_{-}69$.}}
\author {Florin Panaite\\
Institute of Mathematics of the 
Romanian Academy\\ 
PO-Box 1-764, RO-014700 Bucharest, Romania\\
e-mail: Florin.Panaite@imar.ro
}
\date{}
\maketitle

\begin{abstract}
We show that some more results from the literature are particular cases 
of the so-called ''invariance under twisting'' for twisted tensor 
products of algebras, for instance a result of Beattie--Chen--Zhang that 
implies the Blattner-Montgomery duality theorem. 
\end{abstract}
\section*{Introduction}
${\;\;\;\;}$If $A$ and $B$ are (associative unital) algebras and 
$R:B\ot A\rightarrow A\ot B$ is a linear map satisfying certain axioms 
(such an $R$ is called a {\em twisting map}) then $A\ot B$ becomes an 
associative unital algebra with a multiplication defined in terms of 
$R$ and the multiplications of $A$ and $B$; this algebra structure 
on $A\ot B$ is denoted by $A\ot _RB$ and called the 
{\em twisted tensor product} of $A$ and $B$ afforded by $R$ 
(cf. \cite{Cap}, \cite{VanDaele}). 

A very general result about twisted tensor products of algebras was obtained 
in \cite{jlpvo}. It states that, if $A\ot _RB$ is a twisted tensor product of 
algebras and on the vector space $A$ we have one more algebra 
structure denoted by $A'$ and we have also two linear maps 
$\rho , \lambda :A\rightarrow A\ot B$ satisfying a set of conditions, 
then one can define a new map $R':B\ot A'\rightarrow A'\ot B$ by a 
certain formula, this map turns out to be a twisting map and we 
have an algebra isomorphism $A'\ot _{R'}B\simeq A\ot _RB$.  
This result was directly inspired by the invariance under twisting of the 
Hopf smash product (and thus it was called {\em invariance under twisting} 
for twisted tensor products of algebras), but it contains also as 
particular cases a number of independent and previously unrelated 
results from Hopf algebra theory, for instance 
Majid's theorem stating that the Drinfeld double of a quasitriangular Hopf
algebra is isomorphic to an ordinary smash product (cf. \cite{majid}), a
result of Fiore--Steinacker--Wess from \cite{Fiore03a}
concerning a situation where a braided tensor product can be
''unbraided'', and also a result of Fiore from \cite{fiore} concerning a 
situation where a smash product can be ''decoupled''.

The aim of this paper is to show that some more results from the literature 
can be regarded as particular cases of invariance under twisting. Among them is a result 
from \cite{bcz} concerning twistings of comodule algebras 
(which implies the Blattner-Montgomery duality theorem) and a 
generalization (obtained in \cite{lgg}) of  Majid's theorem mentioned before, 
in which quasitriangularity is replaced by a weaker condition, called 
semiquasitriangularity (a concept introduced in \cite{gg1}). 
\section{Preliminaries}\selabel{1}
${\;\;\;\;}$
We work over a commutative field $k$. All algebras, linear spaces
etc. will be over $k$; unadorned $\ot $ means $\ot_k$. By ''algebra'' we 
always mean an associative unital algebra. We will denote by 
$\Delta (h)=h_1\ot h_2$ the comultiplication of a Hopf algebra $H$. 

We recall from \cite{Cap}, \cite{VanDaele} that, given two algebras $A$, $B$ 
and a $k$-linear map $R:B\ot A\rightarrow A\ot B$, with notation 
$R(b\ot a)=a_R\ot b_R$, for $a\in A$, $b\in B$, satisfying the conditions 
$a_R\otimes 1_R=a\otimes 1$, $1_R\otimes b_R=1\otimes b$, 
$(aa')_R\otimes b_R=a_Ra'_r\otimes b_{R_r}$, 
$a_R\otimes (bb')_R=a_{R_r}\otimes b_rb'_R$, 
for all $a, a'\in A$ and $b, b'\in B$ (where $r$ is another 
copy of $R$), if we define on $A\ot B$ a new multiplication, by 
$(a\ot b)(a'\ot b')=aa'_R\ot b_Rb'$, then this multiplication is associative 
with unit $1\ot 1$. In this case, the map $R$ is called 
a {\bf twisting map} between $A$ and $B$ and the new algebra 
structure on $A\ot B$ is denoted by $A\ot _RB$ and called the 
{\bf twisted tensor product} of $A$ and $B$ afforded by $R$. 

\begin{theorem} (\cite{jlpvo}) \label{invtw}
Let $A\ot _RB$ be a twisted tensor product of algebras, and denote
the multiplication of $A$ by $a\ot a'\mapsto aa'$. Assume that on
the vector space $A$ we have one more algebra structure, denoted by
$A'$, with the same unit as $A$ and multiplication denoted by $a\ot
a'\mapsto a*a'$. Assume that we are given two linear maps
$\rho ,\lambda :A\rightarrow A\ot B$, with notation $\rho
(a)=a_{(0)}\ot a_{(1)}$ and $\lambda (a)=a_{[0]}\ot a_{[1]}$, such
that $\rho $ is an algebra map from $A'$ to $A\ot _RB$, $\lambda
(1)=1\ot 1$ and the following relations hold, for all $a, a'\in A$:
\begin{eqnarray}
&&\lambda (aa')=a_{[0]}*(a'_R)_{[0]}\ot (a'_R)_{[1]}(a_{[1]})_R,
\label{4.7} \\
&&a_{(0)_{[0]}}\ot a_{(0)_{[1]}}a_{(1)}=a\ot 1, \label{4.8}\\
&&a_{[0]_{(0)}}\ot a_{[0]_{(1)}}a_{[1]}=a\ot 1. \label{4.9}
\end{eqnarray}
Then the map
$R':B\ot A'\rightarrow A'\ot B$, $R'(b\ot a)=(a_{(0)_R})_{[0]}\ot
(a_{(0)_R})_{[1]}b_Ra_{(1)}$, 
is a twisting map and we have an algebra isomorphism
$A'\ot _{R'}B\simeq A\ot _RB, \;\;a\ot b\mapsto a_{(0)}\ot
a_{(1)}b$.
\end{theorem}

Given an algebra $A$, another algebra structure $A'$ on 
the vector space $A$ (as in Theorem \ref{invtw}) 
may sometimes be obtained by using the following result: 

\begin{proposition}  (\cite{jlpvo}) \label{prop4.1}
Let $A, B$ be two algebras and $R:B\ot A\rightarrow A\ot B$ a linear
map, with notation $R(b\ot a)=a_R\ot b_R$, for all $a\in A$ and
$b\in B$. Assume that we are given two linear maps, $\mu :B\ot
A\rightarrow A$, $\mu (b\ot a)=b\cdot a$, and $\rho :A\rightarrow
A\ot B$, $\rho (a)= a_{(0)}\ot a_{(1)}$, and denote
$a*a':=a_{(0)}(a_{(1)}\cdot a')$, for all $a, a'\in A$. Assume that
the following conditions are satisfied: 
\begin{eqnarray}
&&\rho (1)=1\ot 1,\;\;\; 1\cdot
a=a, \;\;\; a_{(0)}(a_{(1)}\cdot 1)=a, \label{4.1}\\
&&b\cdot (a*a')=a_{(0)_R}(b_Ra_{(1)}\cdot a'), \label{4.2}\\
&&\rho (a*a')=a_{(0)}a'_{(0)_R}\ot a_{(1)_R}a'_{(1)}, \label{4.3}
\end{eqnarray}
for all $a, a'\in A$ and $b\in B$. Then $(A, *, 1)$ is an
associative unital algebra. 
\end{proposition}

\section{The examples}
\setcounter{equation}{0}
\subsection{Twisting comodule algebras}
${\;\;\;\;}$Let $H$ be a finite dimensional Hopf algebra and $A$ a right $H$-comodule 
algebra, with multiplication denoted by $a\ot a'\mapsto aa'$ and comodule 
structure denoted by $A\rightarrow A\ot H$, $a\mapsto a_{<0>}\ot a_{<1>}$. 
Let $\nu :H\rightarrow End (A)$ be a convolution invertible linear map, with convolution 
inverse denoted by $\nu ^{-1}$. For $h\in H$ and $a\in A$, we denote 
$\nu (h)(a)=a\cdot h\in A$. For $a, a'\in A$ we denote 
$a*a'=(a\cdot a'_{<1>})a'_{<0>}\in A$. Assume that, for all $a, a'\in A$ and 
$h\in H$, the following conditions are satisfied:  
\begin{eqnarray}
&&a\cdot 1_H=a, \;\;\; 1_A\cdot h=\varepsilon (h)1_A,   \label{N} \\
&&(a \cdot h_2)_{<0>}\ot (a\cdot h_2)_{<1>}h_1=
a_{<0>}\cdot h_1\ot a_{<1>}h_2,  \label{relalpha} \\
&&(a*a')\cdot h=(a\cdot a'_{<1>}h_2)(a'_{<0>}\cdot h_1). \label{relbeta}
\end{eqnarray}

Then, by \cite{bcz}, Proposition 2.1, $(A, *, 1_A)$  is also a right $H$-comodule algebra (with the same 
$H$-comodule structure as for $A$), denoted in what follows by $A_{\nu }$, and 
moreover $\nu ^{-1}$ satisfies the relations (\ref{relalpha}) and (\ref{relbeta}) for 
$A_{\nu }$, that is, for all $a, a'\in A$ and $h\in H$, we have 
\begin{eqnarray}
&&(\nu ^{-1}(h_2)(a))_{<0>}\ot (\nu ^{-1}(h_2)(a))_{<1>}h_1=
\nu ^{-1}(h_1)(a_{<0>})\ot a_{<1>}h_2, \label{relgamma} \\
&&\nu ^{-1}(h)(aa')=\nu ^{-1}(a'_{<1>}h_2)(a)*\nu ^{-1}(h_1)(a'_{<0>}). \label{reldelta}
\end{eqnarray}
\begin{theorem} (\cite{bcz}) \label{dual}
There exists an algebra isomorphism $A_{\nu }\# H^*\simeq A\# H^*$. 
\end{theorem}

We will prove that Theorem \ref{dual} is a particular case of Theorem \ref{invtw}. 

We take in Theorem \ref{invtw} the algebra $A$ to be the original $H$-comodule 
algebra $A$, the second algebra structure $A'$ on $A$ to be the comodule algebra 
$A_{\nu }$, and $B=H^*$. We consider $A\# H^*$ as the twisted tensor product 
$A\ot _RH^*$, where $R:H^*\ot A\rightarrow A\ot H^*$, $R(\varphi \ot a)=
\varphi _1\cdot a\ot \varphi _2=a_{<0>}\ot \varphi \leftharpoonup a_{<1>}$, 
for all $\varphi \in H^*$ and $a\in A$, where $\leftharpoonup $ is the right 
regular action of $H$ on $H^*$. Define the map $\rho :A_{\nu }\rightarrow 
A\# H^*$, $\rho (a)=\sum _ia\cdot e_i\# e^i:=a_{(0)}\ot a_{(1)}$, where 
$\{e_i\}$ and $\{e^i\}$ are dual bases in $H$ and $H^*$. We will prove 
that $\rho $ is an algebra map. First, by using (\ref{N}), it is easy to see that 
$\rho (1_A)=1_A\# \varepsilon $. We prove that $\rho $ is multiplicative. 
For $a, a'\in A$, we have: 
\begin{eqnarray*}
\rho (a*a')&=&\sum _i(a*a')\cdot e_i\ot e^i\\
&\overset{(\ref{relbeta})}{=}&\sum _i 
(a\cdot a'_{<1>}(e_i)_2)(a'_{<0>}\cdot (e_i)_1)\ot e^i, 
\end{eqnarray*}  
which applied on some $h\in H$ on the second component 
gives $(a\cdot a'_{<1>}h_2)(a'_{<0>}\cdot h_1)$. 
On the other hand, we have 
\begin{eqnarray*}
\rho (a)\rho (a')&=&\sum _{i, j}(a\cdot e_i\# e^i)(a'\cdot e_j\# e^j)\\
&=&\sum _{i, j} (a\cdot e_i)((e^i)_1\cdot (a'\cdot e_j)\# (e^i)_2e^j),
\end{eqnarray*}
which applied on some $h\in H$ on the second component gives \\[2mm]
${\;\;\;\;\;\;\;\;\;\;\;}$
$\sum _i (a\cdot e_i)((e^i)_1 (e^i)_2(h_1)\cdot (a'\cdot h_2))$
\begin{eqnarray*}
&=&\sum _i (a\cdot e_i)(e^i((a'\cdot h_2)_{<1>}h_1)(a'\cdot h_2)_{<0>})\\
&=&(a\cdot (a'\cdot h_2)_{<1>}h_1)(a'\cdot h_2)_{<0>}\\
&\overset{(\ref{relalpha})}{=}&(a\cdot a'_{<1>}h_2)(a'_{<0>}\cdot h_1), 
\end{eqnarray*}
showing that $\rho $ is indeed multiplicative. 

Define now the map $\lambda :A\rightarrow A\ot H^*$, $\lambda (a)=\sum _i 
\nu ^{-1}(e_i)(a)\ot e^i:=a_{[0]}\ot a_{[1]}$. First, it is obvious that 
$\lambda (1_A)=1_A\ot \varepsilon $, because $\nu ^{-1}$ satisfies 
also the condition (\ref{N}). We need to prove now that the relations 
(\ref{4.7}), (\ref{4.8}) and (\ref{4.9}) are satisfied.  It is easy to 
prove (\ref{4.8}) and (\ref{4.9}), because $\nu ^{-1}$ is the 
convolution inverse of $\nu $. We prove now now 
(\ref{4.7}). We have $\lambda (aa')=\sum _i \nu ^{-1}(e_i)(aa')\ot e^i$, 
which applied on some $h\in H$ on the second component gives 
$\nu ^{-1}(h)(aa')$. On the other hand, we have 
\begin{eqnarray*}
a_{[0]}*(a'_R)_{[0]}\ot (a'_R)_{[1]}(a_{[1]})_R
&=&a_{[0]}*(a'_{<0>})_{[0]}\ot (a'_{<0>})_{[1]}(a_{[1]}\leftharpoonup 
a'_{<1>})\\
&=&\sum _{i, j}\nu ^{-1}(e_i)(a)*\nu ^{-1}(e_j)(a'_{<0>})\ot e^j
(e^i\leftharpoonup a'_{<1>}), 
\end{eqnarray*}
which applied on some $h\in H$ on the second component gives 
$\nu ^{-1}(a'_{<1>}h_2)(a)*\nu ^{-1}(h_1)(a'_{<0>})$, and this 
is equal to $\nu ^{-1}(h)(aa')$ because of the relation (\ref{reldelta}). 
Thus, all hypotheses of Theorem \ref{invtw} are fulfilled, so we 
obtain the twisting map $R':H^*\ot A_{\nu }
\rightarrow A_{\nu }\ot H^*$, which looks as follows: 
\begin{eqnarray*}
R'(\varphi \ot a)&=&(a_{(0)_R})_{[0]}\ot (a_{(0)_R})_{[1]}\varphi _Ra_{(1)}\\
&=&a_{(0)_{<0>_{[0]}}}\ot a_{(0)_{<0>_{[1]}}}(\varphi \leftharpoonup 
a_{(0)_{<1>}})a_{(1)}\\
&=&\sum _i (a\cdot e_i)_{<0>_{[0]}}\ot (a\cdot e_i)_{<0>_{[1]}}
(\varphi \leftharpoonup (a\cdot e_i)_{<1>})e^i\\
&=&\sum _{i, j}\nu ^{-1}(e_j)((a\cdot e_i)_{<0>})\ot 
e^j(\varphi \leftharpoonup (a\cdot e_i)_{<1>})e^i,
\end{eqnarray*} 
which applied on some $h\in H$ on the second component gives \\[2mm]
${\;\;\;\;\;}$
$\sum _i \nu ^{-1}(h_1)((a\cdot e_i)_{<0>})\varphi ((a\cdot e_i)_{<1>}
h_2)e^i(h_3)$
\begin{eqnarray*}
&=&\nu ^{-1}(h_1)((a\cdot h_3)_{<0>})\varphi ((a\cdot h_3)_{<1>}h_2)\\
&\overset{(\ref{relalpha})}{=}&\nu ^{-1}(h_1)(a_{<0>}\cdot h_2)
\varphi (a_{<1>}h_3)\\
&=&\nu ^{-1}(h_1)(\nu (h_2)(a_{<0>}))\varphi (a_{<1>}h_3)\\
&=&a_{<0>}\varphi (a_{<1>}h).
\end{eqnarray*}
Thus, we obtained $R'(\varphi \ot a)=a_{<0>}\ot \varphi \leftharpoonup a_{<1>}$, 
for all $\varphi \in H^*$ and $a\in A$, that is $R'=R$ and 
$A_{\nu }\ot _{R'}H^*=A_{\nu }\# H^*$, and so Theorem 
\ref{invtw} provides the algebra isomorphism $A_{\nu }\# H^*\simeq 
A\# H^*$, $a\ot \varphi \mapsto a_{(0)}\ot a_{(1)}\varphi=
\sum _i a\cdot e_i\ot e^i\varphi $, which is exactly Theorem \ref{dual}.
\subsection{External homogenization}
${\;\;\;\;}$Let $H$ be a Hopf algebra and $A$ a right $H$-comodule algebra, with comodule 
structure denoted by $a\mapsto a_{(0)}\ot a_{(1)}$. We also denote 
$a_{(0)}\ot a_{(1)}\ot a_{(2)}=a_{(0)_{(0)}}\ot a_{(0)_{(1)}}\ot a_{(1)}=
a_{(0)}\ot a_{(1)_1}\ot a_{(1)_2}$.  The external homogenization of $A$, introduced 
in \cite{npvo} and denoted by $A[H]$, is an $H$-comodule algebra structure on 
$A\ot H$, with multiplication $(a\ot h)(a'\ot h')=aa'_{(0)}\ot S(a'_{(1)})ha'_{(2)}h'$. 
By \cite {npvo}, $A[H]$ is isomorphic as an algebra to the ordinary tensor product 
$A\ot H$. 

We want to obtain this as a consequence of Theorem \ref{invtw}, 
actually, we will see that the data in Theorem \ref{invtw} lead naturally to the 
multiplication of $A[H]$. Indeed, we will apply Theorem \ref{invtw} to the 
following data: $A$ is the original comodule algebra we started with, $B=H$, $R$ is the 
usual flip between $A$ and $H$, $A'=A$ as an algebra, 
$\rho $ is the comodule structure of $A$ and 
$\lambda :A\rightarrow A\ot H$ is given by $\lambda (a)=a_{(0)}\ot S(a_{(1)}):=
a_{[0]}\ot a_{[1]}$.  It is very easy to see that the hypotheses of Theorem \ref{invtw} 
are fulfilled, so we obtain the twisting map $R':H\ot A\rightarrow A\ot H$ 
given by 
\begin{eqnarray*}
R'(h\ot a)&=&(a_{(0)})_{[0]}\ot (a_{(0)})_{[1]}h a_{(1)}\\
&=&a_{(0)_{(0)}}\ot S(a_{(0)_{(1)}})ha_{(1)}\\
&=&a_{(0)}\ot S(a_{(1)})ha_{(2)}, 
\end{eqnarray*}
and obviously $A\ot _{R'}H=A[H]$. Thus, as a consequence of 
Theorem \ref{invtw}, we obtain the algebra isomorphism from \cite{npvo}:
$A[H]\simeq A\ot H$, 
$a\ot h\mapsto a_{(0)}\ot a_{(1)}h$.
\subsection{Doubles of semiquasitriangular Hopf algebras}
${\;\;\;\;}$Let $H$ be a finite dimensional Hopf algebra and $r\in H\ot H$ an invertible 
element, denoted by $r=r^1\ot r^2$, with inverse $r^{-1}=u^1\ot u^2$. 
Consider the Drinfeld double $D(H)$, which is the tensor product $H^*\ot H$ 
endowed with the multiplication $(\varphi \ot h)(\varphi '\ot h')=
\varphi (h_1\rightharpoonup \varphi ' \leftharpoonup S^{-1}(h_3))\ot h_2h'$, for all 
$h, h'\in H$ and $\varphi , \varphi '\in H^*$, where $\rightharpoonup $ and 
$\leftharpoonup $ are the regular actions of $H$ on $H^*$. 

Define the maps 
\begin{eqnarray*}
&&f:D(H)\rightarrow H^*\ot H, \;\;\; f(\varphi \ot h)=\varphi \leftharpoonup 
S^{-1}(u^1)\ot u^2h, \\
&&g:H^*\ot H\rightarrow D(H), \;\;\;g(\varphi \ot h)=\varphi \leftharpoonup 
S^{-1}(r^1)\ot r^2h.
\end{eqnarray*}
It is obvious that $f$ and $g$ are linear isomorphisms, inverse to each other, so we can 
transfer the algebra structure of $D(H)$ to $H^*\ot H$ via these maps. It is natural 
to ask under what conditions on $r$ this algebra structure on $H^*\ot H$ is a 
twisted tensor product between $H$ and a certain algebra structure on $H^*$. 

We claim that this is the case if $r$ satisfies the following conditions: 
\begin{eqnarray}
&&\Delta (r^1)\ot r^2={\cal R}^1\ot r^1\ot {\cal R}^2r^2, \label{SQT1} \\
&&r^1\ot \Delta (r^2)={\cal R}^1r^1\ot r^2\ot {\cal R}^2, \label{SQT2} \\
&&{\cal R}^1\ot {\cal R}^2_2r^1\ot {\cal R}^2_1r^2=
{\cal R}^1\ot r^1{\cal R}^2_1\ot r^2{\cal R}^2_2, \label{SQT3}
\end{eqnarray}
where ${\cal R}^1\ot {\cal R}^2$ is another copy of $r$. We will obtain this result 
as a consequence of Theorem \ref{invtw}, combined with Proposition \ref{prop4.1}.  Note that the above 
conditions are part of the axioms of a so-called semiquasitriangular structure (cf. \cite{gg1}), 
and that if $r$ satisfies also the other axioms in \cite{gg1} then it was proved in 
\cite {lgg} that $D(H)$ is isomorphic as a Hopf algebra to a Hopf crossed 
product in the sense of \cite{gg2}.  

We take $A=H^*$, with its ordinary algebra structure, $B=H$, and $R:H\ot H^*
\rightarrow H^*\ot H$, $R(h\ot \varphi )=h_1\rightharpoonup \varphi \leftharpoonup 
S^{-1}(h_3)\ot h_2$, hence $A\ot _RB=D(H)$. Then define the maps 
\begin{eqnarray*}
&&\mu :H\ot H^*\rightarrow H^*, \;\;\;\mu (h\ot \varphi )=h\cdot
\varphi:=h_1\rightharpoonup \varphi \leftharpoonup S^{-1}(h_2), \\
&&\rho :H^*\rightarrow H^*\ot H, \;\;\;\rho (\varphi )=\varphi _{(0)}
\ot \varphi _{(1)}:=
\varphi \leftharpoonup S^{-1}(r^1)\ot r^2, \\
&&\lambda :H^*\rightarrow H^*\ot H, \;\;\;\lambda (\varphi )=\varphi _{[0]}
\ot \varphi _{[1]}:=\varphi 
\leftharpoonup S^{-1}(u^1)\ot u^2.
\end{eqnarray*}
The corresponding product $*$ on $H^*$ provided by Propositin \ref{prop4.1} 
is given by
\begin{eqnarray*}
\varphi *\varphi '&=&\varphi _{(0)}(\varphi _{(1)}\cdot \varphi ')\\
&=&(\varphi \leftharpoonup S^{-1}(r^1))(r^2\cdot \varphi ')\\
&=&(\varphi \leftharpoonup S^{-1}(r^1))(r^2_1\rightharpoonup \varphi
' \leftharpoonup S^{-1}(r^2_2)).
\end{eqnarray*}
We need to prove that the relations (\ref{4.1})--(\ref{4.3}) hold. We note first 
that as consequences of (\ref{SQT1}) and (\ref{SQT2}) we obtain 
$\varepsilon (r^1)r^2=r^1\varepsilon (r^2)=1=
\varepsilon (u^1)u^2=u^1\varepsilon (u^2)$, hence we have 
$\rho (\varepsilon )=\lambda (\varepsilon )=\varepsilon \ot 1$ and also 
we obtain immediately $1\cdot \varphi =\varphi $ and 
$\varphi _{(0)}(\varphi _{(1)}\cdot \varepsilon )=\varphi $, 
for all $\varphi \in H^*$, thus (\ref{4.1}) holds. We prove now 
(\ref{4.2}). We compute: 
\begin{eqnarray*}
h\cdot (\varphi * \varphi ')&=&h_1\rightharpoonup (\varphi * \varphi ')
\leftharpoonup S^{-1}(h_2)\\
&=&(h_1\rightharpoonup \varphi \leftharpoonup S^{-1}(h_4r^1))
(h_2r^2_1\rightharpoonup \varphi '\leftharpoonup 
S^{-1}(h_3r^2_2)),
\end{eqnarray*}
\begin{eqnarray*}
\varphi _{(0)_R}(h_R\varphi _{(1)}\cdot \varphi ')&=&
(\varphi \leftharpoonup S^{-1}(r^1))_R(h_Rr^2\cdot \varphi ')\\
&=&(h_1\rightharpoonup \varphi \leftharpoonup S^{-1}(h_3r^1))
(h_2r^2\cdot \varphi ')\\
&=&(h_1\rightharpoonup \varphi \leftharpoonup S^{-1}(h_4r^1))
(h_2r^2_1\rightharpoonup \varphi '\leftharpoonup 
S^{-1}(h_3r^2_2)), \;\;\;q.e.d.
\end{eqnarray*}
In order to prove (\ref{4.3}), we prove first the following relation: 
\begin{eqnarray}
r^1\ot r^2_1\ot r^2_3{\cal R}^1\ot r^2_2{\cal R}^2=
{\cal R}^1_2r^1\ot r^2_1\ot {\cal R}^1_1r^2_2\ot {\cal R}^2. \label{auxil}
\end{eqnarray}
We compute (denoting by $r=\mf R^1\ot \mf R^2=\rho ^1\ot \rho ^2$ two 
more copies of $r$): 
\begin{eqnarray*}
r^1\ot r^2_1\ot r^2_3{\cal R}^1\ot r^2_2{\cal R}^2
&\overset{(\ref{SQT2})}{=}&\mf R^1r^1\ot r^2\ot \mf R^2_2{\cal R}^1
\ot \mf R^2_1{\cal R}^2\\
&\overset{(\ref{SQT3})}{=}&\mf R^1r^1\ot r^2\ot {\cal R}^1\mf R^2_1
\ot {\cal R}^2\mf R^2_2\\
&\overset{(\ref{SQT2})}{=}&\mf R^1\rho ^1r^1\ot r^2\ot {\cal R}^1\rho ^2
\ot {\cal R}^2\mf R^2, 
\end{eqnarray*}
\begin{eqnarray*}
{\cal R}^1_2r^1\ot r^2_1\ot {\cal R}^1_1r^2_2\ot {\cal R}^2
&\overset{(\ref{SQT1})}{=}&
\mf R^1r^1\ot r^2_1\ot {\cal R}^1r^2_2\ot {\cal R}^2\mf R^2\\
&\overset{(\ref{SQT2})}{=}&\mf R^1\rho ^1r^1\ot r^2\ot {\cal R}^1\rho ^2
\ot {\cal R}^2\mf R^2, 
\end{eqnarray*}
and we see that the two terms coincide. Now we prove (\ref{4.3}); we compute:
\begin{eqnarray*}
\rho (\varphi *\varphi ')&=&(\varphi *\varphi ')\leftharpoonup S^{-1}({\cal R}^1)\ot 
{\cal R}^2\\
&=&(\varphi \leftharpoonup S^{-1}({\cal R}^1_2r^1))(r^2_1\rightharpoonup 
\varphi '\leftharpoonup S^{-1}({\cal R}^1_1r^2_2))\ot {\cal R}^2, 
\end{eqnarray*}
\begin{eqnarray*}
\varphi _{(0)}\varphi '_{(0)_R}\ot \varphi _{(1)_R}\varphi '_{(1)}&=&
(\varphi \leftharpoonup S^{-1}(r^1))(\varphi '\leftharpoonup S^{-1}({\cal R}^1))_R
\ot r^2_R{\cal R}^2\\
&=&(\varphi \leftharpoonup S^{-1}(r^1))(r^2_1\rightharpoonup \varphi '\leftharpoonup 
S^{-1}(r^2_3{\cal R}^1))\ot r^2_2{\cal R}^2, 
\end{eqnarray*}
and the two terms are equal because of (\ref{auxil}). 

Thus, we can apply Proposition \ref{prop4.1} and we obtain that 
$(H^*, *, \varepsilon )$ is an associative algebra, which will be denoted in what 
follows by $\underline{H}^*$. 

We will prove now that the hypotheses of  Theorem \ref{invtw} are 
fulfilled, for $A'=\underline{H}^*$. Note first that the relations 
(\ref{4.1}) and (\ref{4.3}) proved before imply that $\rho $ is 
an algebra map from $\underline{H}^*$ to $H^*\ot _RH$. We have 
already seen that $\lambda (\varepsilon )=\varepsilon \ot 1$, so 
we only have to check the relations (\ref{4.7})--(\ref{4.9}). To prove 
(\ref{4.7}), we compute (we denote $r^{-1}=U^1\ot U^2=
\mf U^1\ot \mf U^2$  some more copies of $r^{-1}$): 
\begin{eqnarray*}
\lambda (\varphi \varphi ')&=&(\varphi \varphi ')\leftharpoonup S^{-1}(u^1)
\ot u^2\\
&=&(\varphi \leftharpoonup S^{-1}(u^1_2))(\varphi '\leftharpoonup S^{-1}(u^1_1))
\ot u^2\\
&\overset{(\ref{SQT1})}{=}&(\varphi \leftharpoonup S^{-1}(u^1))
(\varphi '\leftharpoonup S^{-1}(U^1))
\ot u^2U^2, 
\end{eqnarray*}
\begin{eqnarray*}
\varphi _{[0]}*(\varphi '_R)_{[0]}\ot (\varphi '_R)_{[1]}(\varphi _{[1]})_R
&=&(\varphi \leftharpoonup S^{-1}(u^1))*(\varphi '_R\leftharpoonup 
S^{-1}(U^1))\ot U^2u^2_R\\
&=&(\varphi \leftharpoonup S^{-1}(u^1))*
(u^2_1\rightharpoonup \varphi '\leftharpoonup 
S^{-1}(U^1u^2_3))\ot U^2u^2_2\\
&=&(\varphi \leftharpoonup S^{-1}(r^1u^1))
(r^2_1u^2_1\rightharpoonup \varphi '\leftharpoonup 
S^{-1}(r^2_2U^1u^2_3))\ot U^2u^2_2\\
&\overset{(\ref{SQT2})}{=}&(\varphi \leftharpoonup S^{-1}(r^1u^1\mf U^1))
(r^2_1u^2\rightharpoonup \varphi '\leftharpoonup 
S^{-1}(r^2_2U^1\mf U^2_2))\ot U^2\mf U^2_1\\
&\overset{(\ref{SQT3})}{=}&(\varphi \leftharpoonup S^{-1}(r^1u^1\mf U^1))
(r^2_1u^2\rightharpoonup \varphi '\leftharpoonup 
S^{-1}(r^2_2\mf U^2_1U^1))\ot \mf U^2_2U^2\\
&\overset{(\ref{SQT2})}{=}&(\varphi \leftharpoonup S^{-1}(r^1\mf U^1))
(\varphi '\leftharpoonup 
S^{-1}(r^2\mf U^2_1U^1))\ot \mf U^2_2U^2\\
&\overset{(\ref{SQT2})}{=}&(\varphi \leftharpoonup S^{-1}(u^1))
(\varphi '\leftharpoonup 
S^{-1}(U^1))\ot u^2U^2, 
\end{eqnarray*}
and we see that the two terms are equal. The remaining relations 
(\ref{4.8}) and (\ref{4.9}) are very easy to prove and are left 
to the reader. Thus, we can apply Theorem \ref{invtw} and we 
obtain the twisting map $R':H\ot \underline{H}^*\rightarrow \underline{H}^*\ot H$, 
\begin{eqnarray*}
&&R'(h\ot \varphi )=(\varphi _{(0)_R})_{[0]}\ot (\varphi _{(0)_R})_{[1]}
h_R\varphi _{(1)}=h_1\rightharpoonup \varphi \leftharpoonup 
S^{-1}(u^1h_3r^1)\ot u^2h_2r^2, 
\end{eqnarray*}
and the algebra isomorphism $\underline{H}^*\ot _{R'}H\simeq 
H^*\ot _RH=D(H)$, given by
\begin{eqnarray*}
&&\varphi \ot h\mapsto \varphi _{(0)}\ot \varphi _{(1)}h=
\varphi \leftharpoonup S^{-1}(r^1)\ot r^2h,
\end{eqnarray*} 
which is exactly the linear isomorphism $g$ defined before. Thus, we 
have proved that if $r$ satisfies the conditions (\ref{SQT1})--(\ref{SQT3}) 
then $D(H)$ is isomorphic as an algebra to a 
twisted tensor product between $\underline{H}^*$ and $H$. 

\end{document}